%
%
%
\newif\ifmakeref 
\makereffalse    
%
%
%
\baselineskip = 15.5pt plus .5pt minus .5pt   
\hsize = 16truecm
\font\bigbigbf=cmbx10   scaled \magstep1     
\font\sans=cmss10                         
\font\smallcaps = cmcsc10                 
%
\ifx\fonts\cmfonts
\font\ninerm=cmr9
\font\ninei=cmmi9
\font\ninesy=cmsy9
\font\ninebf=cmbx9
\font\ninett=cmtt9
\font\nineit=cmti9
\font\ninesl=cmsl9
\else
\font\ninerm=amr9
\font\ninei=ammi9
\font\ninesy=amsy9
\font\ninebf=ambx9
\font\ninett=amtt9
\font\nineit=amti9
\font\ninesl=amsl9
\fi
\skewchar\ninei='177
\skewchar\ninesy='60
\skewchar\ninett=-1
\newskip\tglue
\def\ninepoint{\def\rm{\fam0\ninerm}
       \textfont0=\ninerm
       \textfont1=\ninei
       \textfont2=\ninesy
       \textfont\itfam=\nineit \def\it{\fam\itfam\nineit}%
       \textfont\slfam=\ninesl \def\sl{\fam\slfam\ninesl}%
       \textfont\ttfam=\ninett \def\tt{\fam\ttfam\ninett}%
       \textfont\bffam=\ninebf \def\bf{\fam\bffam\ninebf}%
       \tt\tglue=.5em plus.25em minus .15em
       \normalbaselineskip=11pt
       \setbox\strutbox=\hbox{\vrule height8pt depth3pt width0pt}%
       \let\sc=\sevenrm \let\big=\ninebig \normalbaselines\rm}%
\input amssym
\def\REFERENCES{recsyn.ref}
\input \REFERENCES
\def\fonts{cmfonts}
\advance\tolerance by 100

%
\def\runningtitlestring{Recursive Synthesis and Foundations}
\newcount\footnotecount
\footnotecount = 0
\def\makeftn[[#1]]{~\hskip-.3em{\ninepoint%
\global\advance\footnotecount by 1
\footnote{\hskip-.3em$^{\number\footnotecount}$}{#1}}}
\newcount\sectioncount
\sectioncount = 0

\def\section#1{\vskip0pt plus .1\vsize
    \penalty-250\vskip0pt plus-.1\vsize\bigskip
    \global\advance\sectioncount by 1
    \centerline{\bf \number\sectioncount. #1}\medskip\message{#1}}
\def\SectionBreak{\vskip0pt plus .1\vsize
    \penalty-250\vskip0pt plus-.1\vsize\null\bigskip}
\def\subsection#1{\vskip0pt plus .01\vsize%
\penalty-250\vskip0pt plus-.01\vsize\bigskip%
\noindent{\bf #1.\ }\message{#1}}

\def\abstract#1{{\ninepoint\bigskip\centerline{\hbox{
       \vbox{\hsize=4.85truein{\noindent ABSTRACT.\enspace{#1}}
            }}}}}

\def\AmsClassMM#1{{\ninepoint\bigskip\centerline{\hbox{
       \vbox{\hsize=5.75truein{\noindent
       AMS 2000 Subject Classification Numbers:\enspace#1}
            }}}}}
\def\th #1 #2: #3\par{\medbreak{\bf#1 #2:
\enspace}{\sl#3\par}\par\medbreak}
\def\co #1 #2: #3\par{\medbreak{\bf#1 #2:
\enspace}{\sl#3\par}\par\medbreak}
\def\le #1 #2: #3\par{\medbreak{\bf #1 #2:
\enspace}{\sl #3\par}\par\medbreak}
\def\rem #1 #2. #3\par{\medbreak{\bf #1 #2.
\enspace}{#3}\par\medbreak}

\def\sqr#1#2{{\vcenter{\hrule height.#2pt
      \hbox{\vrule width.#2pt height#1pt \kern#1pt
         \vrule width.#2pt}
       \hrule height.#2pt}}}

\overfullrule=10pt
\def\boxit#1{\vbox{\hrule \hbox{\vrule \kern2pt
                 \vbox{\kern2pt#1\kern2pt}\kern2pt\vrule}\hrule}}
%
\newdimen\refindent\newdimen\plusindent
\newdimen\refskip\newdimen\tempindent
\newdimen\extraindent
\newcount\refcount
\newwrite\reffile
\def\beginref{\ifmakeref\immediate\openout\reffile = \REFERENCES\else\fi}
\def\endref{\ifmakeref\immediate\closeout\reffile\else\fi}
\def\ifundefined#1{\expandafter\ifx\csname#1\endcsname\relax}
\def\referto[#1]{\ifundefined{#1}[?]\else[\csname#1\endcsname]\fi}
%
%
\refcount=0
\def\ref#1:#2.-#3[#4]{\ninepoint 
\advance\refcount by 1
\setbox0=\hbox{[\number\refcount]}\refindent=\wd0
\plusindent=\refskip\extraindent=\refskip
\advance\plusindent by -\refindent\tempindent=\parindent %
\parindent=0pt\par\hangindent\extraindent %
 [\number\refcount]\hskip\plusindent #1:{\sl#2},#3
\parindent=\tempindent
\ifmakeref
\immediate\write\reffile{\string\expandafter\def\noexpand\csname#4\endcsname%
{\number\refcount}}\else\fi}
\refskip=\parindent
%

\def\makeknown#1#2{\expandafter\gdef
            \csname#2\endcsname{\hbox{\csname#1\endcsname #2}}}
\makeknown{rm}{Cov}
\makeknown{rm}{int}
\makeknown{rm}{Ord}
\makeknown{rm}{Rank}
\makeknown{rm}{rank}
\makeknown{rm}{Root}
\makeknown{rm}{End}
\makeknown{rm}{Level}
\makeknown{rm}{Seg}
\makeknown{rm}{pr}
\makeknown{rm}{top}
\makeknown{rm}{cf}
\makeknown{rm}{loc}
\makeknown{rm}{Hom}
\makeknown{rm}{Pt}
\makeknown{bf}{Loc}
\makeknown{rm}{op}
\makeknown{bf}{Top}
\makeknown{bf}{p}
\makeknown{sans}{COV}

%

\def\specialheadlines{
        \headline={\vbox{\line{
                     \ifnum\pageno<2 \ffolio
                     \else\rightheadline
                     \fi}\bigskip
                     }}}
        \def\ffolio{\hfil}
        \def\rightheadline{\hfil
                    {\smallcaps\runningtitlestring}\hfil
                    \hbox{\rm\folio}}
\nopagenumbers
\specialheadlines
%
%

\def\atiyah{\makeftn[[Our description may be compared
with that in \referto[Atiyaha], which contrasts `geometric thinking' with `algebraic
manipulation'. The promise
of the algebraic machinery - especially since Descartes introduced
his algebraic coordinates - is described as an "Faustian offer":
"...when you pass over into algebraic calculation ... you stop
thinking geometrically ... you stop thinking about the meaning."
(\referto[Atiyaha], p.\ 7).
One sells one's soul in return for obtaining the elements necessary
(in our analysis)
for comparing and manipulating all concepts in a uniform, algebraic
fashion. And here, of course, the `algebra of sets' is the ultimate
reduction of the `geometry' of recursive synthesis.]]}
\null
\vskip1.5truecm
  \centerline{\bigbigbf Recursive Synthesis and the Foundations of Mathematics}
\vskip.75truecm
 \centerline{\sl Aarno Hohti, University of Helsinki}
\vskip.5truecm
\vskip.5truecm
\abstract{This paper presents mathematics as a general science of
          computation in a way different from the tradition. It is
          based on the radical philosophical standpoint according
          to which the content, meaning and justification of experience
          lies in its precise formulation. The requirement on precise,
          formal content discloses the relational structure of
          (mathematical) experience, and gives a new meaning to the
          `ideal' objects beyond concrete forms. The paper also 
          provides a systematic reason why set theory represents 
          an ultimate stage in mathematical technology.}

%
\smallskip
\AmsClassMM{Primary 03A05, secondary 00A30.}
\vskip.5truecm
%

\section{Introduction}

Philosophies of mathematics fall into established trends such as
Platonism, intuitionism, constructivism, modalism,
and others. They are usually studied on the basis
of a special `mathematical experience'. The radical
standpoint adopted here seeks to bring the foundational
study of mathematics to the global philosophical context
of {\it experience in general}. Traditional foundational
frameworks (e.g., set theory, category theory, proof theory)
start from special internal fields within mathematics,
and they attempt at describing the entire, present and future,
mathematical enterprise from their perspective. Critically
speaking, they try to `force' the rest of mathematics into
their mold. The purpose of this paper is to view
mathematical foundations from the standpoint of a
general structure or `logic' of experience.

Our goal is expressed in the need of grounding the
philosophy of mathematics on {\it expressible } and
{\it intercommunicable} forms, rather than
`personal intuition'. We request a critical
evaluation of experience, both `natural' and
`mathematical'. The actual (formal) content of
experience is to be set against the remainder
that cannot be expressed in public, written
forms. Philosophical and even more so, scientific,
considerations should be based on, or `effectively'
connected with, formal experience. A fortiori,
this ought to be expected in mathematics.

The point of departure here lies in a rigorous
{\it articulation} of experience. We say that
experience is {\it bounded} by its forms;
each element of experience
is a form, and as such is determinate and bounded
(not merely finite). While the form is bounded
in itself, it is -- in general -- related to other
forms. Thus, a particular form lies in a
{\it relational field}. The field is not there
as a totality of its elements but rather as
a {\it possibility} or mode of forms to be connected to a
given one. In general, the totality would not be 
defined within the given formal experience.

I am concretely limited with respect to my
capabilities: For example, I may easily add
small numbers, say 3 and 4, without the use
of an algorithmic method of counting. However,
for larger numbers, we introduce new {\it units}
that replace groups with the same number of elements.
Thus, we replace ten elements by a unit, and
repeat this as many times as needed to express
the desired number with nine `digits'. {\it Recursive
synthesis}, simply described, consists in replacing
the many with one, an element to be used in
further syntheses. This new element is {\it related}
to the replaced ones through the very relation
of collecting (synthesis). It both 'hides' and
'opens up into' the collected elements.

Instead of beginning from the relational field of
concrete, formal experience, {\it analytic philosophy}
has based its arguments on a `free' space of
sets and individuals, placed at an ideal separation
from concrete experience. Set theory, quantification
theory and analytic philosophy are mutually related
through their conceptual dependence on `collections
of individuals'. The first two have become the methodology
of analytic philosophy, and they share the same
{\it epistemological notion of experience}.

For recursive synthesis, the principal notion is
that of collecting a multitude under a unit,
i.e., forming a {\it totality}. Set theory
has placed its notion of totality in
`pure epistemology'\makeftn[[Pure epistemology places
its objects outside of experience and description --
those are objects of which we can have {\it knowledge}.]].
The problematic
totalities are the infinite ones. As shown by
various approaches to this field, the
philosophy of mathematics hinges on the
notion of infinity\makeftn[[This does not preclude the 
existence of other such `hinges'.]]. Here we take, however,
our point of departure by asking for the
precise meaning of `finite' objects. We do
not take them for granted.

The simple notion of finitude is {\it transitive}:
A finite collection of finite collections is
again finite. This notion corresponds
to the simple infinity of pure epistemology,
completed in itself. In fact, transitivity is a
form of completion or closure. Transitivity
is the essence of topology in both the
classical sense (Kuratowski closure) and
the Grothendieck topologies, and ties these
together with modal logic (system S4). Transitivity
naturally leads to the completion, as in the
case of positive integers one is led to the
minimal collection which contains all the repeated
successors of the first element. However, the notion of
`all the successors' (in its common, non-critical sense)
is not anymore given in
formal experience but is epistemological. It is
postulated as a well-known object of mathematical
thought, but it is not formally represented
\makeftn[[Axiomatic set theory has to assume the existence
of an infinite set and then
derives the existence of a set containing
the given finite sets as subsets.
But in postulating a {\it set}, it goes
beyond the {\it formally given}, inasmuch as
sets have intuitive `content' at all.]].
On the
other hand, in a purely relational setting,
we are given a form $f$ and a formal definition
describing to which
other forms $f$ is related. The problem is to
derive a necessary (and sufficient) description of
the `totality' of forms, thus related, without going
beyond the given experience. It turns out that
the synthetic `form' of the totalities is found
in their topological `structure'.

What are then the `concrete' forms of mathematical experience?
In this paper, forms are anything that can be given
a finite alphabet and  finite rules of generation;
in fact, a generative or `productive' relation.
For example, we may regard the formulas of a standard
first-order predicate logic in this manner. (In general, 
we should not stipulate that the forms be `generated' - 
it should be enough to check that a formal relation
holds between any two {\it given} forms.)

Each part of mathematics has its own productive rules
for the forms, and its internal relation which make
the forms into a relational field. Thus, the forms
in a propositional calculus are its sentences, generated
in the standard manner from the logical connectives and
propositional symbols. The associated relational field
is obtained from the relation of implication. There is
no need for individuals.

In a relational field, each element is {\it bounded}
by its relational neighbourhood. Similarly,
to be bounded
determines a relationship between the given and
the other, in other words, the boundary is always
between an `interior' and an `exterior'. This is
the natural point of contact with topology and modal logic.
Indeed, modal logic has been re-appearing
in philosophical and logical considerations through
centuries, and recently as a general problem
for the analytic logic that made Quine to discard
it altogether. From our standpoint, the reason
that the {\it ghost of modality} (as coined by
Weyl \referto[Weyla]) keeps re-instating its relevence,
is to be found in the way it is associated with
the structure of relational fields.

Sentential modal logic (i.\ e., without
quantifiers) itself is a part of the theory of relational
fields. By the work of Kripke , we know that for
each set with a relation, there is a corresponding
`modality' that defines a modal logic, and vice versa
(the elements of the set being `possible worlds').
In our setting, the modality is the
relation which describes what forms are `possible'
for a given one. The link between modal logic
(specifically, the so-called system S4) and topology has been
understood since the work of MacKinsey and Tarski
in the 1930's\makeftn[[A related link between
sentential intuitionistic logic and modal logic (S4)
was pointed out by G\"odel \referto[Godela] in 1933.]].

In fact, {\it covering systems} of forms (or propositions)
--- or more generally partially ordered sets ---
satisfying natural weak axioms were considered by
Fourman and Grayson \referto[FourmanGrayson] to define
{\it formal topologies}\makeftn[[For an exposition of
this theory, see also Sigstam \referto[Sigstama].]].
A covering relation holds between an element
of a partially ordered set $P$ and a subset $A\subset P$.
Of the four axioms (equivalent to those
of a Grothendieck topology) in \referto[FourmanGrayson],
three are rather obvious\makeftn[[i.\ e., 1) a subset
covers each of its elements 2) if $a \leq b$, then
$\{b\}$ covers $a$, and 3) the `meet' of two sets
covering an element still covers the element.]],
but the fourth --- postulating the {\it transivity}
of the concept of covering --- is essential from
the standpoint of topology. Transitivity corresponds
to the specific axiom $MMp\rightarrow Mp$ of S4 in
modal logic, as well as to the idempotency $C^2 = C$
of the Kuratowski closure operator in topology.

The axioms for topology are defined for `parts' of a
topological `space'. In our situation, we
consider {\it chains} in relational fields.
In the most general setting, a chain is a pair $(f,R)$
consisting of a form $f$ and a relation connecting
$f$ to other forms. The relation $R$ is not necessarily
that of the underlying relational field; in general,
it imposes additional conditions on the forms in the chain
and therefore is {\it stronger} than the underlying
relation. The category-theoretical counterpart of a
chain is a projective (or `inverse') system (often denoted as
$(x_i\gets x_j)$; the category theory associates
with each such system the projective limit
$\lim_\gets x_i$. Similarly, we associate a {\it projective
foundation} with a chain, determined topologically
from the members of the chain.
As such, the projective foundation (intuitively: limit)
of the chain is nothing but its `place'. It has no
content in itself other than its topological
structure. Traditionally, the projective foundation
is the {\it set} of such chains (for example, the real
numbers are given as a set of (equivalence classes of)
Cauchy sequences). The topological meaning of the 
projective foundation is the main reason that separates our 
views from 
set theory as well analytic philosophy. 

In our work, however, sets are replaced by
{\it concepts} (`generic
elements', `parts') represented by forms in
the given relational system. Even in classical,
set-theoretical mathematics, sets are frequently
defined as collections of elements that satisfy a
formula. Usually in such situations, the set does not
represent anything over and above the formula itself.
However, often set-theoretical considerations
go beyond formal experience and involve sets
which have no formal definition (sufficient to
identify them) but are justified on
epistemological grounds.

\section{A Critique of `Pure Epistemology'}

A criticism of a tradition of rational thinking that
had surpassed the limits of possible experience was already
presented by Kant. The famous antinomies (whether well-founded
or not) are reminiscent of the dilemmas of set theory.
We are not thinking here of Russel's paradox, but
rather results of independence, e.g., that of the
Continuum Hypothesis. By going from finite formulas
to objects of epistemology\makeftn[[i.e., in this
situation to objects (sets) which are stated to be there,
but which have no {\it explicit} formal definition,
such as the elements of the power set of an infinite
set.]] -- which are not controlled by these formulas --
it is not surprising that we find antithetical
statements between which we cannot decide.

Kant wanted to limit our arguments about experience to
possible experience. But this domain is confused, in the
sense that finitude as such is `transitive'. Husserl's
first work ({\sl Philosophie der Arithmetik}) was
characterized by the idea of bounded experience, shown
by the distinction between direct and `figurative'
collecting of elements into unities. This early
attempt became considered `psychologistic', and Husserl
later provided other (more philosophically `correct')
pathways toward the same goal: Phenomenological reduction,
ideation, and the relationality of the {\it Lifeworld}.
But in our opinion, the phenomenological approach still lacks (Husserl
threw the same criticism against Kant) a {\it rational basis}.
It was not based on well-defined, intercommunicable
formal experience -- which necessarily is concretely bounded
and relational. (Nevertheless, it shares with us the
understanding that a basis for philosophy is to be found in
experience itself.)

The radical standpoint of {\it formal experience} should be contrasted
with that of epistemology, which considers {\it objects of knowledge},
not those of experience. Thus, one may treat infinite sets as objects
of knowledge, even if we cannot experience them as such. This does not
mean that mathematicians do not possess a special skill of `seeing'
or visualizing infinite collections. However, the concept of set
is changed when we move from finite, elementwise articulated
collections to infinite ones. In the latter case it is the defining
formula of the set that precedes the elements -- these are present
only in the sense of an epistemological `completion'. The problem
at hand is not about these completions as such, but rather pertains
to their nature.

In accepting the Power Set Axiom, we leave the domain within
which every element (`individual') is separately and explicitly
formally definable. Whatever our countable language may be,
there will be elements of ${\cal P}(N)$ which cannot be generated or
defined by any algorithm let alone a finite finite formula
of the language. This is the {\it decisive step} to the purely
`epistemological' sets.

Clearly, standard mathematical practice requires that there be
a {\it `power object'}. But is it necessary to take this
power object to be a set? This is a crucial question for set theory.
It hinges on the power set axiom.
A successful introduction of set theory
requires that mathematical objects such as $\Bbb R$, the real numbers,
be represented as sets. Higher cardinalities are obtained via the
power set axiom, yielding ever-higher cardinals. 

However,
even the axiomatic set theory only really provides countable mathematical
objects. The L\"owenheim-Skolem `paradox' has been explained away by saying
that an object is uncountable because it has no mapping onto
$\omega$.
The result of L\"owenheim and Skolem says: Every consistent (countable) theory has a
countable model. {\it But what do we mean by an uncountable
model?} How can we even prove that it is {\it genuinely} uncountable?
(In other words, not merely in the sense that the model lacks
an enumerative mapping.)
One needs a foundation on which to establish such genuine
uncountability in the first place. But how does one prove that this
foundation itself is genuinely uncountable?

The only way in which set theory can provide uncountable
sets is by the method of proof, i.\ e., by showing that
for a particular set $S$, the cardinality of $S$ is greater
than that of $\Bbb N$. On the other hand, it is possible
to construct {\it countable models} in which such a proof
is valid. Indeed, a {\it critical} examination {\it shows}
that there are {\it no other} models to begin with.
In the beginning, we have to start with a countable language,
and the very first uncountable set is obtained from a
{\it statement} proving its uncountability. However, this kind of
uncountability is {\it virtual}: We have no other
way of establishing that the first-mentioned model is uncountable
than the theorem itself. But in fact the theorem only states
that the model does not have an enumerating mapping. Let us
assume that the first uncountable set has been obtained by
applying the Power Set construction to a countable set.
This construction adds
a countable collection of canonical elements
of the power set plus the assurance that the new set is
uncountable, i.\ e., there is no
contradicting enumerating mapping.

For the `na\"\i ve' mathematician, this standard resolution of the
L\"owenheim-Skolem problem by a deficiency
is not sufficient. Na\"\i ve thinking
considers all sets in their final (completed) sense. If a model
of the real numbers, say, is countable,
then for the na\"\i ve thinker, it is so in an absolute sense.
In the axiomatic setting, the model is uncountable because it
lacks an enumerating mapping. This non-existence `confirms'
Cantor's proof of the non-countability of the real numbers.
But there is, of course, such a function, necessarily {\it outside}
of the model. But for the na\"\i ve thinker, this function should
be {\it as good as any} --- there are no models but only the absolute set
of real numbers. If you can somehow enumerate a set, inside the model
or not, then the set is enumerable. As the conclusion, we obtain that
axiomatic set theory, rigorously speaking,
cannot provide a full account of {\it the set}
of real numbers, it merely gives partial, countable `views' which
are not always mutually consistent. No model of the real numbers cannot be proved
to be genuinely uncountable, and a merely virtually uncountable
model is --- for the na\"\i ve thinker --- already countable.
So the common mathematician has to reach beyond axiomatic
set theory to the na\"\i ve one, where paradoxes such as Russell's
are waiting.

The L\"owenheim--Skolem dilemma makes us suspicious of infinite sets
{\it in general}. Indeed, if the only {\it known} genuinely infinite sets are the
countable ones, we should investigate the reason for this
rather special situation. Countable sets can be
'presented' through their enumeration.
What makes such an object different from a
mere algorithm of this construction is the possibility of
treating the object as a {\it totality}, as a unity, instead
of considering it element by element, produced from the
algorithm. (And the totality is distinct from the form which
is the algorithm.) The convenience in epistemology is to
consider the totality a set, amenable to the constructions
used with finite sets (such as the power set), again producing
other sets (of higher and higher `cardinalities'). What is
needed here --- according to our radical standpoint --- is not
to discard these infinite unities as such. Rather, we must
critically examine the {\it concrete} mathematical `experience'
(or rather, `presentation')
of such collections, and the relationship of the collecting
unity to the elements collected.


Why do we consider set theory to have failed as a foundation for mathematics? Let us
consider the following points.

{

%

\item{1)} From the beginning, it has been beset either by
  paradoxes or undecidable questions;
  \smallskip

\item{2)} It adds an unnecessary purely epistemological or
  `metaphysical' layer of transfinite sets to the formalism
  of mathematics;
  \smallskip

\item{3)} With the Skolem-L\"owenheim problem, as explained above,
  even the motivation of advanced set theory (as that of transfinite numbers)
  is in doubt.
  There is no way to prove that there are genuinely uncountable sets.
  \smallskip

\item{4)} To even prove that there are countably infinite {\it sets}, one
  still needs an axiom. This is a purely
  epistemological assumption\makeftn[[This axiom might be called
  metaphysical, or `ontological' (Russell), but as its object is
  the mathematical realm of abstract knowledge, we call it
  `epistemological'.]]. Also, it is the second `leg' of transfinite 
  set theory. 
  \smallskip

\item{5)} Set theory postulates objects (as simple as real numbers) of
  which one cannot say anything descriptive at all.

}

\noindent
But does not {\it every} real number have a decimal representation?
However, given an arbitrary (abstract) real number $x$, the
decimal series $x_1,x_2,x_3,\ldots$ can only be given if there is
an algorithm (rule) that generates the sequence. Hence, there are
uncountably many (in fact, `almost all') real numbers destined to
remain unknown forever, unless our mathematical language somehow can
be extended to encompass an uncountable variety of finite forms.
It is mysterious how such an extension could be achieved. It is,
however, even more intellectually troubling
that this residual set of real numbers may not, after all,
be genuinely uncountable, but nevertheless
resists our attempts at describing the elements!

By the Skolem-L\"owenheim result, our collection of real numbers
may (assumed to) be based on a countable model. But what could
such a {\it countable} region of inaccessibility be? It is easy
to understand why in uncountable regions, some elements remain
indistinguishable. On the other hand, for a countable set, we
should be able to identify the elements, one after another, in
an enumerative process. Nevertheless, these elements remain hidden.
The reason, of course, is that we don't know where to look for
them, we don't know where to start and how to continue. The
countability of the model underlying this uncountable remainder
is beyond the model itself. This conflicts with the na\"ive,
absolute point of view which takes each set in a
fundamental, `final' way: If there is a method, perhaps in a larger
model, to enumerate a set, then it is countable. The conlusion
is we cannot speak of the set of real numbers as such,
but only of {\it models} of real
numbers, which should be disappointing to those who had expected
set theory fulfil the role of a foundation for mathematics.

As a remedy to the critical remarks 1) -- 5) above (and
against the recommendation of \referto[Putnama]) we seek
in recursive synthesis a foundation which provides mathematics
a place in a unity of sciences rather than a
position separated from the rest. But it will turn out that
this foundation cannot be {\it mathematical}. Indeed, if
mathematics is a study of general forms of recursive
synthesis, then it itself has no general form, and hence no
foundation within itself.

\section{The Uniformization of Recursive Synthesis}

The production of forms provides as such no ground
for the logical movement from an argument to another;
it merely adds to the multiplicity of increasingly
special situations. Moreover, the finite
structure of the logical space does not
solve the {\it recursively synthetic problem of
movement} beyond the limits of the bounded
reason. This problem of {\it combination}
cannot be solved by means of a multitude of natural forms;
to become effective, the collecting process
of synthesis must be grounded on a simple
basis which reduces the {\it general field} of
such problems. 

The problem of overcoming
boundedness, in whatever field, is to connect
two points $A, B$ that have been separately presented
within the limits of our cognition. But this
relating requires a {\it comparison} of $A$
with $B$, it calls for a {\it ground of movement}
from the former to the latter. The initial
limitedness is shown in the mere {\it juxtaposition}
of the two which is the simple synthesis
placing $A$ and $B$ side by side without a ground.
This juxtaposition is `groundless' as such, shown in the meaning
of the synthesis that opens up to this
opposited juxtaposition and hence to the
emptiness of its connection. The synthetic
relation can be grounded only if it opens up to
a {\it common foundation} of both parts, in other
words only if we can find a {\it uniform grounding}
of both $A$ and $B$.

The idea of finding a common basis in geometrical
arguments was known in the Greek tradition
in its conception of {\it analysis and synthesis}.
To connect $A$ and $B$, one must bring them onto a common ground
through their analysis to simple parts, shared or
comparable, and then proceed by synthesis along
the path thus paved, indeed by reversing the
process of analysis. Uniformization begins when
the fixed units of such divisions are isolated
and named, and continues recursively when the results of basic operations between
them become replaced by known entities of a
scale.
It is precisely
in the introduction of the {\it symbolic algebra} that
recursive synthesis shows itself: a group of
expressions is replaced by a new `symbol'. This
is an object to be handled in the same way as the
others; it is a member of a
{\it homogeneous}
field of objects. Such uniformities had been
known in the
Greek mathematics\makeftn[[The Greeks did not accept even rational numbers as such;
instead of a `real' fraction they had proportions
of integers. Numbers other than integers and proportions
were given by geometric figures, by concrete
synthetic forms in which the various number
constituents were joined. Each number had a definite
form, obtained from simpler forms, and the questions of
their `totality' were meaningless in that context.
They were natural,
total forms ordered with respect to their recursive
collecting of their parts ({\sans `methexis'}). This hierarchy
has a simple expression in the geometric scale of dimensions.
The reduction of
this recursive synthesis through symbolic algebra
and the subsequent simplification by means of `coordinates' constitutes
a paradigmatic example of the process of uniformization.
(As regards
the {\it `geometric algebra'} of the Greeks,
the reader is referred to B.\ L.\ van der Waerden:
{\sl Science Awakening}, P.\ Noordhooff, 1954.)]]
as basic elements
corresponding to various dimensions: line, `square' ({\sans dynamis}),
`cube' ({\sans kybos}) and their later higher-dimensional extensions
such as {\sans dynamodynamis} etc. Each uniformity
forms a {\it species} within which one computes
by addition, because the various objects of
any one such field are multiples of the same
basic unit, in the sense that lengths, areas and
volumes each admit of addition. This understanding
became systematic in Vi\`ete's {\it logistice speciosa},
i.\ e., calculus of species, a general algebra
of indeterminate entities to be contrasted with
the {\it logistice numerosa} of common numerical
computations. Indeed, in the {\sl Analytic Art}
Vi\`ete enounced his {\it law of homogeneity},
the idea of calculus with magnitudes within
uniform domains (genera) for which there has been
set ``a series or ladder [{\it series seu scala}] ... of magnitudes
ascending or descending by their own nature from genus
to genus''\makeftn[[{\sl Opera Mathematica}, Georg
Olms Verlag, 1970, p.\ 1. Translated in
Klein \referto[Kleina] p.\ 322 (see also
pp.\ 172--173).]]. At each level of the scale, the law of homogeneity
permits addition and subtraction, while division and
multiplication (as Vi\`ete shows in Chapter IV of
the {\sl Analytical Art}) operate {\it between}
these recursive steps.

The {\it general foundation of calculus} is the uniformity
of its objects with respect to each other. Thus,
the result $5+2 = 3+4$ is possible with respect to a
uniform domain of integers based on a simple unit.
Recursive synthesis is made possible by the condition that
every result of finite collecting is again an element
of the domain, directly comparable and combinable with any other
element. In particular, algorithmic
computation based on a decimal system becomes possible
by the equivalence of all collections of ten units and
so on to higher and higher synthetic levels of collecting.
The efficiency of the
positional system lies in the uniformity of
scaling: all the positions are mutually similar;
the method of position moves the uniformity of the original
unit (one) to that of the different levels
(powers of ten).


Uniformization is essential for a mathematics of recursive
synthesis. The early uses of calculus replaced
the `real objects' by uniform units which were mutually
equivalent from the standpoint of calculation. The integers
are uniform units of collecting, and they abstract from the
non-essential features of the collections. At this level,
they are effectively merely {\it representational}.
In order to make collecting and computing effective,
a uniform scale of recursive levels --- powers of 10 ---
was introduced. By using this scale, the addition of any two
integers could be reduced to repeated additions of two
digits (one had to know the `addition table' by heart).
An essential element in uniformized recursive synthesis,
and perhaps its simplest form, is the use of {\it products}.
In fact, even the simplest collecting of elements under
a fixed element can be characterized as a cone product
of the collection and the collecting element. More
important is the case of a product between two non-trivial
factors. Products --- or conversely division --- are uniform
reformulations of syntheses, and uniform scales of products
(powers) make a uniformized recursive synthesis possible.

Finding a uniform foundation of basic elements (`units')
is merely the first (but necessary) step in bringing a
variety of given objects into framework within which comparison
and calculus is possible. This was essentially
understood by Descartes (discussed in more detail
in the next section), and as he says, it was already
exemplified in the method
of analysis and synthesis in geometry. But in order to make
the ensuing mathematics of synthesis effective, and not
simply remaining content with haphazard combinations of units,
the recursive synthesis should be uniformized, too.
A mathematical field typically has a particular kind
of synthetic relation to be uniformized. For example,
for Vi\'ete it was synthesis by multiplication.

If the greatest achievements in mediaeval mathematics were
characterized by algebra --- abstract and recursive use
of multiplication and addition --- the Renaissance introduced
a method of multiplicative simplification which became the
essential mathematical technique for hundreds of years.
Differential calculus starts from {\it linearization},
i.\ e., replacing a curve or a function by a straight line
or a linear function. One considers
expressions of the form
$ f(x + h) \; = \;f(x) + a h + \epsilon(h) h, $
where $a$ is the `derivative' of $f$ at $x$ and $\epsilon h$
is the difference between the linear function
$g: x\mapsto a (x+h)$ and $f$. Ideally speaking, $f$
is replaced by a linear function in infinitesimally small
neighbourhoods of $x$. We may call this the
{\it first linearization} of $f$.

Recursively, we next target the difference, which already
contains a multiplicative term. We `linearize' the non-linear
part of the difference, and obtain an expression of the form
$ f(x+h) \;=\; f(x) + a h + bh^2 + \epsilon_2(h) h^2,$
the {\it second linearization} of $f$.
In general, we derive the rule for the Taylor development
of $f$ at $x$ expressing the function as a series in the
successive powers $x^n$ of $x$, which now plays the role
of a uniform scale.

Classical predicate calculus itself can be
described in this manner!
In the situation of a general relational field, still
considering the classical case in which we have
a finite predicate language and a model, synthesis amounts
to expressing the relational neighbourhoods of the
elements in the model. The distributive normal
forms and constituents of Hintikka (\referto[Hintikkaa])
disclose the relational and synthetic meaning hidden in
the formulas of predicate logic, which traditionally
treats relations from the standpoint of set theory.
Constituents of degree zero
are simply exhaustive combinations of atomic formulas
describing which atomic properties hold for a
given $k$-tuple of elements. For degree 1, we describe
what kind of elements `exist' (in our setting: are related to)
for the given $k$-tuple, for example,
$$C^{(1)}_*(y,x) = B_j(x,y) \& \Pi_s\pm\exists z C^{(0)}_s(x,y,z),$$
is a constituent of degree 1 for the pair $(x,y)$, where
$B_j(x,y)$ is a combination of atomic formulas and their
negations that hold for the pair. Each constituent is a
conjunct ('product') of constituents of a lower level.
The above formula leads
to the general rule for defining constituents of any degree.
The constituents $C^k$ represent {\it uniformized units of
description} in the sense that {\it any} formula $\varphi(x)$ of
the underlying language is equivalent to a unique
disjunction of such constituents $C^d(x)$, where $d$
is the `quantifier depth' of $\varphi$. What is more,
for a given element $a$ of the model, there is a
unique chain (series)
$\ldots\gets C^d_i(x)\gets C^{d+1}_j(x)\gets C^{d+2}_k(x)\gets\ldots$
of constituents which represents the maximal description
of $a$ by using formulas with one free variable.

\section{The Idea of Simplification in Descartes and Boole}

An answer to the problem of a complete uniformization of both
planar and spatial forms is 
obtained by a simplification
known as the {\it introduction of coordinates}.
The projection of a geometric figure onto a line
is a {\it valuation} which `simplifies' the {\it geometric
point} to a number and thus makes the points
comparable with each other. The figures
considered are not anymore restricted to the
complexes of natural forms; the permitted
figures are now complexes of the basic
valuations (coordinate symbols) and `constants'.
The recursive hierarchy of forms -- for which
scales were used by Vi\`ete -- has now
been transformed into a homogeneous field
of calculus in which individual operations
are simply `coordinatewise arithmetic operations'.

While coordinate representations had been used
earlier, the uniformization of geometry by
adding the simplifying coordinate valuations
to symbolic algebra is connected with Descartes.
However, the case of geometry is merely a particular
example; the {\sl Regulae}\makeftn[[{\sl Regulae
ad Directionem Ingenii}, in {\sl Oeuvres}, Chapter X.
We have used the translation by Marion: {\sl Regles
Utiles et Clares pour la Direction de
l'Esprit en la Recherche de la Verite},
Martinus Nijhoff, 1977.]] shows a more general
understanding of recursive synthesis and its uniformization
as the basis of human knowledge in general. What is
at stake here is the long-disputed deeper meaning
of his {\it `mathesis universalis'} as the basis of the
true method. With respect to
this method, the sciences of Arithmetic and Geometry are
``nothing but ripened fruits'' (373) from the principles
of the deeper mathesis that is nothing like the
common mathematics, {\it `mathematica vulgaris'} (374).

For Descartes,the required true universal
mathesis is not merely ``a certain general science
that explains all that can be investigated concerning
order and measure'' (378); what is required is the
{\it systematic uniformization} of that recursivity which is
needed for the overcoming of our boundedness. This
was not made clear in the full generality even in
Descartes. But he understood that in order to make the movement
of thought ``continuous and nowhere interrupted'' (Rule VII)
the parts must be arranged in series with simple relations
the connectivity of which is immediate and clear. Uniformization
consists in finding the units
 which provide the common ground
of the things considered, the unit is ``the basis and
foundation of all the relations'' (462).
The
method of simplification is included in Rule XVII: The
problems must be transformed until they are reduced
to finding ``certain magnitudes'' (459); indeed, the
relations between things are to be {\it reduced to those between lines}
``because I could find nothing more easily pictured
to my imagination'' (441)\makeftn[[But even the case of
measurement which as such requires mediation
is reducible to order (452) for which the relation of
terms is a direct one (451), simple recursion
in which a term points at the objects immediately
closed in.]].

One needs a
systematic method of {\it comparing} concepts, down
to the `last elements' without remaining uncertainty.
The ultimate method of comparison is to reduce concepts
in a way similar to Descartes' coordinization in geometry,
replacing the real values by `truth values'.

The solution to the problem of conceptual simplification
by means of valuations was to be found by Boole.
His work is nowadays associated with the notion of
`Boolean algebra', a system in which certain rules
crystallized by Boole are valid. However, in order to
understand its essential relation to the uniformization
of synthesis, we must consider the original formulation
of Boole's thoughts in terms of `elective symbols'
which correspond to so-called characteristic functions,
valuations from the `universe' of discourse to the
simple `algebra' $\{0,1\}$.

The truth-functional method of Post and Wittgenstein reduces
the propositional forms to finitely many distributions of these simple
elements, and the `content' of the propositions in those forms
is left out of consideration. Thus, from the {\it point of view
of propositional truth} the form is fully given as the synthesis of
those distributions, as for example the form $p\to q$ is analyzed
into the truth distribution $1,0,1,1$ over the four cases
$(0,0),(1,0),(0,1),(1,1)$.
The synthesis which it provides
for the propositional form, say $p\to q$, from the four binary
valuations, is in the simplest view a synthesis through
a mere union, discrete collection. It provides a method
of testing the validity of the propositional form by going through all
the cases, one by one.
The truth-functional method
is in full accordance with the effective mathesis of the
Cartesian method, the universal technique of uniformizing
recursive synthesis that analyzes
with respect to a simple unit
and puts together results layer by layer in uniform steps.
Propositional truth-functional
calculus carries with itself the ambiguity of both ends, of computation
and of logic in the traditional sense, as if it were at the same time
a {\it method}, a {\it `clavis universalis'} into the
classical forms of logic. But what it effectively and concretely
does is a manipulation of binary assignments; to each {\it binary
input} there corresponds a {\it binary output}. No wonder that Boole was
accused by some of his contemporaries of deducing correct results
by unfounded methods! The pre-Boolean world was still
unaccustomed to the thought that logical forms might be
faithfully handled by transforming them (in terms of `elective functions')
into the simpler realm of zeros and ones, although the method
is internally consistent. But these critics were right insofar
that what really is handled here is the binary `image' of the
traditional logical forms and nothing beyond it; with this step the
classical logic became effectively replaced by the new binary logic.

In general, Leibniz's logic is still a `concept logic' based
on the containment of concepts: The calculus of {\it `continentibus
et contentis'} understands the statement {\it $A$ est $B$} as the
containment of $A$ by $B$ instead of the modern interpretation as
`predication' by $B$, $B(A)$. This is properly speaking neither
extensional nor simply intensional but a combination of both.
Indeed, not only is the `concept' $A$  in $B$, but according to the
extensional interpretation every individual (i.\ e.,
complete concept) `in $A$' belongs to $B$. But this means
{\it all possible individuals}\makeftn[[as pointed out
by W.\ Lenzen \referto[Lenzena], esp.\ pp.\ 168--170.
See also his {\sl Das System der Leibnizschen Logik}, de Gruyter,
Berlin, 1990.]], not merely the `real' ones. Thus, $A$ refers
to the conceptual ground `consisting' of all possible extensions
of $A$. By a theorem of Leibniz\makeftn[[in
{\sl Opuscules et fragments in\'edits de Leibniz},
      edited by L.\ Couturat, Paris, 1903,
      (reprinted by Georg Olms, Hildesheim, 1961), p.\ 260.]]
`$A$ is $B$' is equivalent to the statement `if $X$ is in
$A$, then $X$ is in $B$', in other words the ground of $A$
is a ground of $B$. Quantification was still directly related
to its older sense as `supposition'. The ground is not extensional
in the sense of a set or a class, but rather an `extensive'
variable, and for its conceptual quantification, the statement
{\it some $A$ is $B$} is then given as the identity
$AY = B$, where $Y$ is an indefinite concept compatible with $A$.
Thus, even quantification remains within the field of
conceptual and {\it formal composition}.

The entire Leibnizian project, with its emphasis on the
conceptual understanding of logical forms, with its understanding of
calculus as given by rules of transformation, and with its ideal
of a {\sl Characteristica Universalis} reflecting the true structure
of thoughts and concepts, resists the recursive uniformization
of Descartes. For this very reason, his {\sl Analysis Situs} and more
generally the ideal of
{\sl Characteristica Geometrica} could not advance to an `effective'
level, because by employing such primitive concepts as point and surface
and their intrinsic relations, it could not have reached the
facility that only a deeper analysis into simple elements can bring
before us. It is due to this reason that Grassman's {\sl `Ausdehnungslehre'}
provided a workable geometric calculus because it presents
a construction of its objects based in his `exterior' product.
But Grassman's calculus was not a geometric topic in the sense
of Leibniz; it was not designed as a calculus of `native' geometric
forms and their congruences. Only at the end of the 19th century did
such a calculus become effectively realized through the {\it division}
of those original forms into polyhedra of `simplices', combinatorially
uniform elements of analysis, and through the replacement of their
mutual relations by their indices of `incidence': a number 
which is either $1$, $0$, or $-1$. In this way, the
{\it geometric form} became `valuated' into a matrix of indices
of incidence along the same lines as in the method of `truth tables'
with respect to {\it propositional forms}. It is on this basis
that a systematic and even an automated process can assign
`topological invariants' to the original form.

\section{Uniformization in set theory}

We have considered the uniformization of recursive synthesis
in several particular cases.
In counting, the natural numbers
provided a general foundation for collecting finite sets.
According to Klein \referto[Kleina],
the Greek conception of number in its
deepest sense is that of {\it ordering} a set of elements
by keeping them together in a synthetic unity, shown
in the figurate numbers of the Pythagoreans.
In abstract counting, the order is inessential,
what matters is the finite {\it set} as the
basis of the cardinality. General set theory both
continues counting into the transfinite and at the same time
seeks to provide a uniform synthetic foundation for all
mathematical objects.

The method of uniformization obtains its most
universal and radical expression in the theory of
sets that provides the ultimate grounding by
final elements, mutually equivalent and indivisible
points. The objects of this theory are sets or `classes'
composed of these atomic individuals, the synthesis
of the composition being the simple collecting of the elements,
as in the pure form of finite synthesis. But what is important
here is that set theory does not stop with finite collections;
it extends the principle of gathering a set over all multitudes
of individuals, and postulates in this process the existence of
at least one non-finite set, an object with `infinitely many'
distinguishable elements. In order to proceed at all beyond merely
being a theory of finite sets -- which would by no means reduce
it to a triviality -- it has to assume the existence of a set
that essentially is the `set' of natural numbers, as the basis
of its induction and its generation of `transfinite' ordinal and
cardinal numbers.

The uniformity of the set-theoretical synthesis from the
foundation of natural numbers is expressed in the
`cumulative hierarchy'. It represents
the `iterative concept' of sets. In the hierarchy, sets are
constructed out of the `empty' set by
the repeated use of set-theoretical operations. This
depends on the first  and already existing ground of natural numbers,
without which it could not carry its iteration
beyond finite stages. (More precisely, one
may use any infinite set as a foundation, but
the existence of an infinite set implies that of natural
numbers.) We may as well proceed by applying
the power set operation to the foundation.
With the cumulative hierarchy, the recursive
synthesis of mathematical objects is then uniformized
by the foundation and the power set operation
iteratively applied to the foundation. Let us
notice that the `usual' mathematical objects
are already obtained in the first few levels
of the hierarchy. In the same way as the standard
decimal notation provides a uniformized synthetic
representation of natural numbers, and power series
gives one for the analytic functions, the cumulative
hierarchy provides a framework for {\it all} mathematical
objects.

Uniformization is the {\it ultimate justification} of set theory.
We have discussed above the role of uniformization
as the `technology' of recursive synthesis, and interpreted
Descartes' {\sl Regulae} in this light as the
universal recipe for systematically solving any problem\atiyah.
Therefore, set theory is not merely a useful approach
to the foundations of mathematics, but even
stipulated by the logic of technology, given
our recursively synthetic conception of experience.
Synthesis from a common foundation, following a
uniform method, makes all objects thus represented
mutually comparable. But as with another
example of uniformization --- the gene technology
which provides such a foundation for biology ---
the question is {\it what other} objects
(`monsters') might become produced alongside
with those originally desired? The set-theoretical
synthesis produces (via the power set operation)
sets with no definition; indeed it stipulates
the existence of objects which will forever
remain `in the dark side'. Depending on the
model, even `virtual sets' (Cohen reals, etc.)
can be found.

Perhaps {\it ethical guidelines} should be
applied here, too. They would delimit our powers of
set-theoretical construction, as in the predicative
analysis, eliminating the
bogus syntheses that come with an unfettered freedom.
On the other hand, even the question of foundation
(natural numbers, infinite set, etc.) is problematic.
To provide a complete foundation of units for all
mathematical objects (or at the minimum for the
objects of classical analysis), set theory has to
go beyond all finite levels to a totality which
simultaneously realizes all these levels at once.

\section{The `utopia' of the set-theoretical infinite}

The existence of the foundation by the simultaneous being-there
of natural numbers is absolutely necessary for the development
of transfinite numbers. Indeed, although natural numbers can be
generated from the empty set $\emptyset$ by
an iterated formation of sets, so that the numbers $0,1,2,\ldots$
are given as the sets
$\emptyset,\{\emptyset\},\{\emptyset,\{\emptyset\}\}$ and so on,
there is no infinite generation without the prior existence
of an {\it already} infinite set $\omega$, the ground
and supply of all the further development in set theory.
Therefore, insofar as the otherwise merely finite theory is relatively
free of philosophical problems, it is {\it this assumption
of the first infinity} that requires our attention.

As a simple synthesis,
$\omega$ is not given as a projective ground (foundation),
because the latter depends on things beyond set theory
and is already a {\it complex} synthesis involving the notions of
fundamental neighbourhood and rule (Cf.\ the next section of this paper.)

It is a {\it simple} synthesis of its elements, but the synthesis
is not concretely (formally) given -- it is not a place for this
synthesis (not being localized at anything concrete). Therefore,
it is `placeless'; we call it {\it utopian} ({\sans topos}, `place').


In the name of scientific honesty, one should in the context of set
theory do away with all interpretations of infinity, as was
precisely postulated in Hilbert's formalism of signs.
While Cantor's theory purported to ground all mathematical
objects -- and therewith the vast realm of transfinite numbers --
once and for all by setting a ground of a self-subsistent infinite,
the Hilbertian pure play of finite symbols would eliminate the need of the
resulting inaccurate understanding of the infinite by
leaving it {\it uninterpreted}. As the meaning of the infinite
in mathematics was not yet ``inexhaustibly ({\it restlos})
explained''\makeftn[[Hilbert \referto[Hilberta], p.\ 161.]], one should replace --
following the achievement of Weierstra\ss\ in the
case of mathematical analysis -- the uses of the infinite with
``finite processes that perform the same task'' (162). The
``final explication'' of the ``essence of the infinite''
had become necessary for the ``honour of human understanding''
(163). However, the explication of the infinite by the ``contentless'
use of finite operational symbolism through which its role
is reduced to ``that  of an idea'' (190) does not do
anything to relieve us from the need to assume the principle of infinity
and in fact rather leaves it in the darkness. The result of
this symbolism is again {\it as if} the problem of the infinite
had been resolved -- by looking away from it -- but it has merely
been relegated from mathematics qua mathematics to the
{\it practicing subject}. It is the practical author and the reader
of mathematical works who has to provide the understanding
of that pure symbolism and of setting it in concrete operational
movements, most importantly in the sense of inductive operations
which implicitly prepostulate the infinite in the form
of the `and so on' etc. In fact, what Weierstra\ss\ obtained
for differential and integral calculus was not reproduced
by Hilbert, because the quintessential topo-logical content
of the Weierstra\ss ian `reduction' is neglected
by Hilbert. The so-called $(\epsilon,\delta)$-definition
of continuity and limits produces this topo-logical content, while in
Hilbert the infinite is reduced to the `ideality' of the
pure symbol, which lacks the dynamical character of the Weierstra\ss ian
version of analysis. Again, the ideal element is
`as if' it were infinite, although all that is concretely given is a
finite symbol.

On the other hand, the assumption of an infinite set beyond concrete
experience as a `pure' synthesis might not lead to any
{\it mathematical harm} inasmuch as even the wildest imagination
must be accompanied by a sound mathematical practice. Nevertheless,
it leads to questions that cannot be decided on the basis of concrete
experience, in the same way as the assumption of a {\it divinity
beyond our solid experience} leads to problems that must remain
unresolved. The Renaissance scientist could add the existence of a God
beyond human powers, and without falling into contradictions
with everyday experience could assume various mutually inconsistent
features of the postulated God. The set-theorist has likewise been
left with questions that cannot be solved on the basis of
set-theoretical evidence alone, the most notable of these
problems being the well-known question of the `cardinality'
of the continuum.

By placing the primary infinite object beyond concrete experience,
set theory realizes the ultimate ideal of uniformization by entering
a {\it utopia}, a {\it place which is nowhere}; it is not localized
as a projective ground. The set-theoretical
infinite set is a {\it concrete beyondness} vis-\`a-vis
every possible experience.
In its
original form in Cantor, this utopian thinking is not restricted
to a mere empty {\it displacement} of the ground of the
infinite, but extends to an ideal of `transfinite numbers',
a well-ordered linear hierarchy of everincreasing ordinal numbers
which -- with their rigorously defined arithmetical
rules -- extend the realm of finite ordinals. For Cantor the pure
infinite cannot be expressed in mathematics. However, it was in the
spectrum of the actually (and not merely potentially) infinite
ordinal numbers that mathematics approached the `total' infinite.
The first transfinite number to follow natural numbers is the
ordinal corresponding to $\omega$, i.\ e., the totality
of the natural numbers. The ordinals to follow this
first transfinite one are $\omega+1$, $\omega+2$, ... and so on to ever larger
transfinites. Cantor sought to demonstrate that
{\it infinite numbers} as such are not an {\it inconsistent utopia}
-- as had been thought previously by many -- but
an {\it actual fact} by {\it representing them in concreto}
in the form of the transfinite sequence of concrete
ordinal numbers.

However, the very problem of constructing `transfinite' numbers,
i.\ e., numbers which are `larger' than the ordinary natural
numbers, presents in itself no special difficulty. Any object
distinct from numbers, e.\ g.\ this table or that book
can be declared `transfinite' as soon as as we relate it
in a suitable way to the natural numbers. But it would be
needless to think of such a tangible object as being
`beyond' or `behind' the whole sequence of numbers, in the same
sense that the first transfinite ordinal is waiting there, in its
utopia, to be met after running through the classical number sequence.
Indeed, this concrete object would be {\it with us}, together with
the first integers. And nevertheless it would be in some definite
sense `beyond' the whole collection of natural numbers
by simply not being one of them. This table or this book
is here; I do not have to run through the sequence of
numbers in order to catch a glimpse of it. Therefore,
my placing it after the sequence is an {\it imaginative emplacement}
that has nothing to do with purely mathematical
considerations: the actual relation of this concrete object
with the numbers is merely a relation {\it alongside with}
other relations and the image of `beyondness' is
my imaginative contribution, an ideal of `coming
after'\makeftn[[The original motivation for Cantor's
theory seems to have arisen from infinite iterations of certain operations
for sets of real numbers. Those examples are also sources
of the original `pictures' of ordinal numbers. In a simplified
form, we may take a converging sequence such as
$1, 1/2, 1/3,\ldots$ with the limit $0$ to represent a
natural `model' for the ordinal number $\omega$. After
infinitely many eliminations of the first term of the
sequence, i.\ e., after removing $1, 1/2, 1/3$ etc.\ one
after another, we are still left with the limit term.
In this sense, it comes `after' all the diminishing positive
terms of the sequence. But in a different sense
it is present there in the very same sense
as the other terms that are concretely given. However,
{\it as such} this limit term is nothing but a real number; it becomes
meaningful as the limit or the first term beyond the series
only in the context of that series. In other words, this
situation is present only when the total collection with the
limit is considered. It is the `open' series that is the
focal point here; the limit $0$ is merely one of its
{\it finite determinations} as the end point of the
region of openness.]].

The ordinal number $\omega + 1$ is given
{\it explicitly} and not merely implicitly as e.\ g.\ the number
$2^{73}$ of which we do not immediately have a `full'
presentation. We do not need an infinite development of the
explicit expression $\omega + 1$ in the same way  as
$2^{73}$ requires a finite development (to express it as a
decimal number); indeed, the expression `$\omega + 1$'
is not only the most precise formula for this number;
it is concretely given in its explicitness
here and now. As in the case of this table or some other
`natural' object, $\omega + 1$ is {\it by the side}
along with the positive integers and not behind them; it is
a `number' which can be represented as a pair
$<1,2,3,\ldots;1>$, in which the difference between
`coming after' and `lying beside' is one of imagination,
an addendum that corresponds to the utopia of the infinitely
far.

The representation $<1,2,3,\ldots;1>$ does not introduce a `new'
infinity beyond and above the first infinity $\omega$; indeed, it
merely shows the old infinity with a {\it new form}
having two ends rather than one. In the same way, the representation
of $2\omega$ as $<1,2,3,\ldots;1,2,3,\ldots>$ only involves the
same original openness and a new {\it entirely finite form}
given by the brackets $<;>$ and so on to more and more complex forms. In
fact,
each `countable' ordinal may be represented in a form of a `tree'
with only finite `chains', obtained recursively as a finite form molding
the first infinity already present in $\omega$.
In those forms, there is no more `infinity' than in the first;
what separates them is the complexity of their finite form.
The interval $[0,1]$ is another `formed' infinity in which the
topo-logical structure is given by the rule of its subdivisions.
In the same
way as the move from $I^1$ to the second power $I^2$ does not
add a new infinity but merely {\it complicates} the division
of the same original openness, so does the move from $\omega$
to $\omega+\omega$ just add a new finite dimension.

Constructivism preserves the talk of sets,
but declares that construction is the
criterion of existence, thereby cutting the utopian elements
of Cantor's set theory down to the experience with concrete
content. In Brouwer and
Weyl the topo-logical structure of openness of the continuum
is present in the notions of `free becoming' and
in the representation of real numbers as `choice sequences'
that are in the state of free prolongement in the same way as
any decimal representation of a non-rational real number
is only given up to a finite initial segment together
with the rule of its further expansion.
Gounding is given by `rules', the collection
of which is never complete but always under creation and incrementation.


But the application of the rules of the constructivist and
the very understanding of their idea is left to the {\it human
subject}. By failing to thematize the openness of the generation
of its objects, it has pushed the infinite -- in the way of
Hilbert's symbolism -- away from the mathematical thematic to the
practicing mathematician. It has to presuppose the `first'
infinite of an infinite series as the very principle of its rules,
and it sufficiently describes it by the rule of simple becoming
from $n$ to $n+1$. But while this may be considered satisfactory for the
natural numbers, passing by the topo-logical content shows itself
in the poverty of its description of more geometric objects
such as intervals and circles etc. By not taking the model
of its generative principle from geometry (as Frege during his
last stage) but rather from the discrete generation of
sequences, the `almost invisible' topo-logical structure
of the natural number series easily became ignored in its
{\it philosophy of the mere rule}.

There is a way to do justice to the classical
mathematical intuition of objects such as lines and circles, without
taking away anything of their content, preserving their
concrete presence without advancing to an unreachable
`beyond' of the Cantorian set theory. While the traditional finitism
or constructivism may rigorously exhibit only a partial realization
of the continuum, and the classical set theory enters -- instead
of the concrete infinity of the continuum -- the path of abstraction
beyond experience, the examination of our concrete intuition of
the continuum and other geometric objects discloses
their topo-logical structure in which the infinite openness
(divisibility) is given as a concrete content, at once, and
not as a mere possibility of unlimited division. The
interval $[0,1]$ as the projective limit is {\it not a set of individual
elements} (as in Cantor) nor lacks any real numbers (as in Brouwer),
but it is the projective foundation of the relation of betweenness (or
that of inclusion between subintervals). Instead of amending the
serial philosophy with openness ('choice sequences')\makeftn[[An
individual series is still always connected with a
rule which selects the consecutive successors for the given
relation.]] -- a
difficult marriage between determinism and freedom -- we replace
successor rules with relations. The series are now objects derived
from the relation.
We may still say that there
{\it are} $2^{\omega}$
real numbers in the sense that they form a projective foundation
similar to that of $(2^n\gets 2^{n+1})$. However, at this conceptual level,
we may {\it not} say that some real numbers exist because they have
recursive developments and others do not.

\section{Some Main Notions for a New Foundation}

In this section, we briefly examine notions
involved in developing a `new foundation' based
on the above philosophical considerations. What follows
will merely indicate the path to be taken
in that work. A detailed development of basic
mathematics on the foundations of recursive synthesis
will be presented (attempted) in another article.

\subsection{A. Extension Relations}
Instead of sets, we consider {\it formal concepts}. In the 
basic case, these are given by a {\it form} $f$ and a 
formal relation $R$ that defines extensions to the form.
Traditionally, extension has been separated from the form:
the extension of $F$ has been said to consist of individuals
$x$ such that $F(x)$ (i.e., $x$ `satisfies' $F$). However,
in our situation there are {\it no} individuals distinct
from forms. The `extension' of a form consists of 
forms related to the former through the given extension 
relation. 

For each relation $R$, and for each $f$, there is the
relational neighbourhood $R[f]$ of $g$ with respect
to $R$. However, it is not a set. Inasmuch as $R$ is
formal, $R[f]$ is the formal {\it concept} of an
object $g$ such that $fRg$. For example, let $S$
be the successor relation of the natural numbers.
Then $S[1]$ is the concept of the successor of 1,
which is realized by a unique number. On the other
hand, let $R$ be the inclusion relation between
open rational intervals of the form $]r,q[$.
Thus, $]r,q[R]s,t[$ if, and only if,
$r < s < t < q$. The neighbourhood of $]r,q[$
with respect to $R$ is the concept of an open subinterval
of $]r,q[$. Again, `interval' does not directly imply a
set; here a rational interval is a primitive
form composed of 5 symbols. It is important to notice 
that the forms discussed here are the mathematical 
objects constructed (or rather `formulated') within
a formal language. 

For us, a formal language is a finite alphabet equipped with
a productive relation (which distinguishes well-formed
expressions within the strings generated from the
alphabet).
A simple {\it formal system} ${\cal F}$ is a formal language
together with an internal extension
relation between the forms. To consider a particularly
simple example, the formal system of decimal series
of real numbers in $[0,1[$ consists of the
alphabet $<0,1,2,\ldots,9>$, the
production relation $R$ for which $fRg$ iff $f$
is of the form $0.\alpha_1\cdots\alpha_n$,
and $g$ is of the $f\alpha$, where $\alpha$
is a symbol of the alphabet. The internal
extension relation is the same as the
production relation.

The extension relation gives a rule which, however, is
not computable in general. In the case of integers, it is
the effective rule $n\gets n+1$, as in the
case of recursive series. On the other hand,
the extension relation for real numbers does not yield
an effective functional rule, but nevertheless is given by
a well-defined and effectively decidable
form: $]r_1,r_2[R]s_1,s_2[$ iff
$r_1 < s_1 < s_2 < r_2$.

\subsection{A note on Russell's paradox} 
In order to attempt deriving Russell's paradox in our setting,
we could consider a binary relation $R$ and form the concept
$F$ defined via $Fx\equiv \neg xRx$. Here $F$ is 
is a unary relation (property) of forms. To get the paradox, 
we should have $Fg \equiv gRF$ for forms $g$; this implies
$FRF \equiv \neg FRF$ as in the set-theoretical case. This, 
of course, can be done, witnessed by the set-theoretical 
case and its membership relation ($\in$). However, the 
paradox requires that $R$ be applicable to $F$, while $F$ 
(as a form) has been constructed using the symbol $R$. 
Our requirement that all forms (and hence relations) be 
constructed from the given finite alphabet helps to avoid 
this type of circularity.

Set theory forces a single relation where there should be 
distinct ones. It diagonalizes forth Russell's paradox
by equating (in the above example) the {\it relation} $R$ and 
the {\it property} of satisfying $\neg xRx$. 

Circularity in general is {\it a} motivation behind the study
of non-wellfounded sets, including cycles such as 
$\Omega = \{\Omega\}$\makeftn[[See, in particular, Barwise and 
Moss \referto[Barwisa].]]. However, already the formulation of 
the {\it Axiom of Antifoundation} (every graph has a decoration
by sets) indicates the naturality of graphs and general 
relations v\'\i s-a-v\'\i s the membership relation. In a 
relational formal system, sets are not needed.

\subsection{B. Chains}
We often denote extension relations by
arrows $\gets$.
Given such a relation $\gets$, we may consider objects $h$
such that there is $g$ with $f\gets g$ and $g\gets h$. More
generally, we may consider sequences
$f_1,\ldots,f_n$ such that
$f\gets f_1$, $\ldots$, $f_{n-1}\gets f_n$.
Such objects $f_n$ are called transitively related to
$f$. The {\it transitive completion} of $\gets$,
denoted by $\gets^*$, is again defined by giving
its extension relation. Thus, we define
$$ \gets^n \quad\Leftarrow\quad \gets^{n+1} $$
by the rule: if $f\gets^n g$ and $g\gets h$, then
$f\gets^{n+1} h$. The new relation is the
path relation derived from $\gets$.

Alternatively, we may define the concept of
path. Let $<>$ denote the empty path. The
concept of $\gets$-path is defined by the
relation $R$ for which $fRg$ iff $f$ is of the
form $<f,f'>$, where $f\gets f'$. Thus, a
$\gets$-path (based at $<>$) is of the form
$$ <<<\ldots<f_1,f_2>,\ldots,f_n>,$$
which may be simplified to the form
$<f_1,\ldots,f_n>$.

Relations of a formal system ${\cal F}$ can be ordered
by setting $R' < R$ if $fR'g$ implies $fRg$. An $R$-chain
is simply a relation $R' < R$. We may call a chain $R'$
{\it projective} if 
$fR'g$,$fR'h$ implies $gR'k$ and
$hR'k$ for some $k$.
Next, we define a notion which derives its motivation from
combinatorial topology. We call a chain $R'$ {\it closed}
if $xRy$ implies $yRz$ if there is such a $z$, and 
similarly for $x$: $uRx$ is there is such a $u$.
The concept of a {\it closed} $R$-chain leads to that
of the {\it projective foundation} of ${\cal F}$,
denoted here by $\hat {\cal F}$.

\subsection{C. Projective Foundation}
We agree with \referto[Resnika]\makeftn[[esp.\ p.\ 188. 
Here, however, mathematical objects are not primarily
positions in patterns but rather are posited for the 
sake of recursive synthesis.]] that mathematical objects
are posited. Indeed, recursive synthesis proceeds 
by collecting multitudes under unifying elements. Given 
a list $f_1,\ldots,f_n$ of forms, we may add a new 
form $g$ and a relation $R$ such that $gRf_i$ for 
each i. We call $c = <g,R>$ a {\sl cone} of the given 
forms relative to $R$. The relational neighbourhood $R[g]$ is then the 
formal concept of the given list. Note that for 
a given $R,f$, the concept $<f,R>$ already is the cone 
of $R[f]$. (We also denote by $<g,R>$ the sub-chain
of $R$ {\it starting} at $g$.) 

On the other hand, multitudes are also given 
by means of chains. Given a chain $c = <f,R>$,
the cone $\hat c$ of $c$ collects the chain 
through a relation $\hat R$ that relates $\hat c$
to each finite part $f,f_1,\ldots,f_n$ of the chain.
However, let us note (following a typical mathematical 
procedure of `completions') that it is natural to 
choose for $\hat c$ the chain itself. The chain 
is related to each of its finite (initial) parts  
through forgetting the rest of the chain (projection).
We call $\hat c$ the projective limit or foundation of 
the chain. 

In general, the chain relation $R$ allows for 
multiple `successors' of forms in the chain. Hence,
we have a multitude
of chains as distinct specifications of the given concept.
The projective limit $\hat c$ 
enables us to speak about the `totality' of 
such chains. While we can have distinct closed 
$R$-chains $c_1$,$c_2$, they are both attached to $\hat c$.

For example, in the case of the natural numbers
${\cal N} = <0,S>$ (where $S$ denotes the successor
relation), the projective foundation $\hat {\cal N}$
`consists' of a unique closed chain (the relation $S$
itself). In the case of dyadic rational intervals
$]i/2^k, (i+1)/2^k[$ contained
in $[-1,1]$ -- denote this system by
${\cal Q})$ -- the
closed chains correspond to real numbers: For
a dyadic rational interval $]r,s[$, the closed chains
through the interval form a subconcept which
corresponds to the interval $[r,s]$ in the
real line $[-1,1]$ which is the projective foundation $\hat {\cal Q}$.


The `elements' of the concept $\hat R$ should be
associated with minimal closed $R$-chains  $R'$ as principles or
{\it rules of choice}: Such an element $R'$ implies,
for any $x$, exactly one $y$ such that $xR'y$ provided
there is such a $y$ for $R$. Sometimes such principles
can be found, e.g., in the case of `constructive'
real numbers such as $\sqrt 2$ which give an explicit
projective foundation. However, in general there is no
element of this kind. In set-theoretical thinking,
it is customary to assume the set of all such elements,
or, the set of the associated series. For us, the
projective foundation does not consist of `points'
(although we may distinguish specific, definable
points). We may only indicate its topological
structure by means of the fundamental neighbourhoods
$R'[y]$ of extensions from $y$ (more precisely: the
concept of $(R',y)$-extensions.

Projective foundation is the formal notion
of `infinite recursive synthesis'.
Given a chain $c = <f,R>$, the projective foundation
$p = \hat c$ (the `cone' of the chain\makeftn[[In
category theory, projective limits may be considered
examples of
`cones'. For recursive synthesis, a cone is a
paradigm example of collecting a multitude under a unit.]])
has {\it covers} determined by the form $f$ and the
relation $R$ which describes how covers are refined.
First of all, the trivial cover corresponds to
$f$ itself, and a cover of $f$ is composed of
all `parts' corresponding to $f'$ such that
$fRf'$, and so on. From the viewpoint of classical set-theoretical
intuition, such a part
is the `collection' of subchains starting from $f'$,
and minimal closed subchains may be called `points'.
The connection with classical topology is obtained
from the topology of the `space of points' defined
by these covers, but here it is enough to
consider the {\it concept of cover} derived
directly from the given relation. For a
point $p$, the {\it fundamental neighbourhoods}
are precisely the parts corresponding to the forms
of the chain $p$ itself.

The above definition of cover by means of a refinement
relation is directly associated with the set-theoretical
notion of covers.

As an example, consider a relational system where
the forms are finite decimal series and the relation
connects a form $(x_1,\ldots,x_n)$ with any extending
form $(x_1,\ldots,x_n,x_{n+1})$. The associated projective
foundation then `consists' of all infinite
decimal series, but should not be thought of as
a set. It is an object determined by its topology,
for which a canonical open `part' is a part consisting
of infinite series extending a given finite one,
say having $k$ elements.
A canonical cover consists of such open parts, defined
for a given finite $k$. The essential point here
is that while there are many covers, they are all
covers of the {\it same} object.


\def\extensions{
Covers are to be treated as concepts in the same way as
other sets (in the classical sense). However, it is possible
to consider covers by means of {\it conditions} instead of
extension relations. A natural example is given by
the cover of the real line consisting of all intervals
$]x,y[$ such that $|x^2-y^2|<\epsilon$. Here one has an
instance of the {\it splitting of concepts}: In the set theory,
any object is a set. This uniformization is a main advantage
of the set theory. However, when we leave the set theory, we
may experience the splitting of formerly unified concepts.

\subsection{D. Specification}
Functions (or `maps') are here defined between two concepts. A function
$A\to B$ and more generally a relation between
$A,B$ is a formal concept $R$ such that $xRy$ for some pairs
$(x,y)$, where $x$ (resp.\ $y$) is a {\it specification} of
$A$ (resp.\ $B$). Recall that every concept has its relation
of extension (which may be vacuous). Any concept related
to the former via this extension is called its specification.
Note that specification as such does not make a difference
between the element and part of a concept (in the classical
sense). Also note that for any extension relation, the
transitive completion of the relation gives the associated
transitive specification, as in the case of natural numbers.

Let us consider the specification of real numbers. First, we
define the
concept of a real number as the concept of a decimal sequence.
On the other hand, the concept of the real line, is the
projective limit $\hat {\cal R}$ of the concept of a strictly decreasing
sequence of rational intervals. (For all purposes considered
here, it is sufficient to define the corresponding extension
relation as follows: $]r',s'[$ is a strict successor (subinterval) of
$]r,s[$ if $r < r' < s' < s$ and $s' - r' < 2^{-1}(s-r)$.)

We derive the general concept of a sequence from the trivial relation
$R_*$ which relates any two forms. We represent it as the
triple $(*,R_*,\to)$, where $*$ denotes the trivial form
(`any-form'). Then any given sequence
(concept) is a specification of this universal sequence-concept:
Any relation (resp.\ form) is a natural specification of the trivial
relation (resp.\ form).

The classical `elements' of the projective limit
$\hat {\cal R}$ are the decreasing sequences of the above kind.
Here, we have the concept of an ${\cal R}$-sequence, and the
elements are its specifications. What is the specification relation
in this case? Any particular sequence must be given by a
{\it rule} that identifies, for each member of the sequence,
the successor with respect to the extension (${\cal R}$). Such a
{\it selection} is a subrelation of ${\cal R}$, minimal in the
sense that each element has {\it exactly} one successor.

} 

\extensions

\subsection{E. Conclusion}
Set theory was described above as a natural
uniformization based on the set-theoretical element.
Traditionally, set-theoretical mathematics
builds its objects from a given foundation of
basic elements. The recursive synthesis of relational
fields can be considered in two opposite
directions: One which `closes' in the collected
elements under a synthetic element and another
in which the given element is `opened' up to
those elements. In the latter direction, the
ground is not pre-established but rather
topologically given as a projective foundation.

The work presented here may not be the only means
for obtaining the goal set in the introduction.
The program of obtaining and developing classical mathematics in a
predicative framework, from a countable foundation, started
by Weyl and continued by Feferman et al., using more and more
flexible primitives, may ultimately meet our
goals in a common end result: Finitely and rigorously
determined, and at the same adequate, mathematics.
\bigskip
\SectionBreak
\section{References}
\bigskip
{
\beginref

\ref Atiyah, M: Mathematics in the 20th Century.- Bull.\ London Math.\
     Soc.\ 34 (2002), pp.\ 1--15. (Also in {\sl Amer.\ Math.\ Monthly},
     107, (2001), pp.\ 654--666.) [Atiyaha]

\ref Barwise, J., and L.\ Moss: Vicious Circles.- CSLI Lecture Notes 60,
     CSLI Publications, Stanford, 1996. [Barwisa]

\ref Feferman, S: Infinity in Mathematics: Is Cantor Necessary?.-
     Philosophical Topics 17(2), 1989, pp.\ 23--45 (reprinted in
     {\sl In the light of logic}, Oxford University Press, New York, 1998).
     [Fefermanna]

\ref Fourman, M., and R.\ Grayson: Formal spaces.- The L.\ E.\ J.\
     Brouwer Centenary Symposium, A.\ S.\ Troelstra and D.\ van Dalen
     (eds.), North-Holland, 1982, pp.\ 107--122. [FourmanGrayson]

\ref G\"odel, G: Eine Interpretation des intuitionistischen
     Aussagenkalk\"uls.- Ergebnisse eines mathematischen Kolloquiums 4,
     1933, pp.\ 34-38 (reprinted in {\sl Kurt G\"odel, Collected
     Works} (S.\ Feferman et al., eds.), Vol 1, Oxford University Press,
     1986, pp.\ 301--302). [Godela]

\ref Hilbert, D:  \"Uber das Unendliche.- Math.\ Annalen, 95, 1926,
     pp.\ 161--190. [Hilberta]

\ref Hintikka, J: Distributive Normal Forms in the
      Calculus of Predicates.- Acta Philosophica Fennica, Fasc.\ VI,
      Helsinki, 1953. [Hintikkaa]

\ref Klein, J: Greek Mathematical Throught and the
     Origin of Algebra.- MIT Press, 1968. [Kleina]

\ref Lenzen, W: Concepts and Predicates. Leibniz's Challenge
     to Modern Logic.- The Leibniz Renaissance, Biblioteca di Storia
     della Scienza, 28, 1989, pp.\ 153--172. [Lenzena]

\ref MacLane, S., and I.\ Moerdijk: Sheaves in Geometry and Logic.-
     Springer-Verlag, 1992. [MacLanea]

\ref Putnam, H: Mathematics without Foundations.- Journal of Philosophy
     64:1, 1967, pp.\ 5--22 (reprinted in {\sl The Philosophy of
     Mathematics} (W.\ D.\ Hart, ed.), Oxford University Press, 1996,
     pp.\ 168--184). [Putnama]

\ref Resnik, M.D: Mathematics as a Science of Patterns.- Oxford 
     University Press, 1997. [Resnika]  

\ref Sigstam, I: Formal spaces and their effective presentations.-
     Arch.\ Math.\ Logic 34:4, 1995 pp.\ 211 -- 246 [Sigstama]

\ref Weyl, H: The Ghost of Modality.- Philosophical Essays in Memory
     of Edmund Husserl, Cambridge (Mass.), 1940, pp.\ 278--303 (reprinted
     in {\sl Gesammelte Abhandlungen}, vol.\ 3, (ed.\ K.\
     Chandrasekharan), Springer-Verlag, 1968, pp.\ 684--709). [Weyla]

\endref
}

\vskip2cm
\noindent
{\bf The address:}
\bigskip
{\ninepoint\obeylines
\hskip\refskip University of Helsinki
\hskip\refskip Department of Mathematics
\hskip\refskip Yliopistonkatu 5
\hskip\refskip SF--00100 HELSINKI
\hskip\refskip FINLAND
}

\end